\magnification=1200  
\input amstex  
\documentstyle{amsppt}  
\hoffset=.2true in  
\hsize= 6 true in   
\vsize= 9 true in  
\vcorrection{0.2in}   
\NoRunningHeads   
 

\def\ts{\thinspace}%
\def\NN{{\Bbb N}}%
\def\RR{{\Bbb R}}%
\def\SS{{\Bbb S}}%
\def\TT{{\Bbb T}}%
\def\ZZ{{\Bbb Z}}%
\def\QQ{{\Bbb Q}}%
\def\eps{\varepsilon}%
\def\"#1{{\accent"7F #1\penalty10000\hskip 0pt plus 0pt}} 
  
\scrollmode  
 

\topmatter 
\title Asymptotical flatness and cone structure at infinity
\endtitle

\author Anton Petrunin\footnote"$^\bigstar$"
{This work was done while we enjoyed
our stay at the MPI Leipzig
and finished during a joint week at the
IHES.\hfill{$\,$}} and Wilderich Tuschmann$^\bigstar$
\vskip4mm 
\endauthor

\date{}
\enddate

\abstract
We investigate asymptotically flat manifolds with cone structure
at infinity.
We show that any such manifold $M$ has a finite number of ends,
and we classify
(except for the case dim $M=4$, 
where it remains open if one of the theoretically 
possible cones can actually arise)
for simply connected ends all possible cones at infinity.
This result yields in particular a complete classification of 
asymptotically flat manifolds with nonnegative curvature:
The universal covering of an asymptotically flat $m$-manifold 
with nonnegative sectional curvature is isometric to $\RR^{m-2}\times S$,
where $S$ is an asymptotically flat surface.
\endabstract 
\endtopmatter         
\document 
\baselineskip=.25in


\head 0. Introduction 
\endhead

\noindent 
Let $(M,g)$ be a complete noncompact Riemannian manifold.
Choose a point $p\in M$ and set 
$$A(M)=\limsup_{|px|\to\infty} \ \{|K_x|{\cdot}|px|^2\},$$ 
where $|K_x|$ denotes the maximal 
absolute value of the sectional curvatures at the point $x\in M$. 
 
One easily checks that $A(M)$ 
does not depend on the choice of the reference point $p$,
so that the quantity $A(M)$ yields a nice geometric invariant of $M$
which is, in particular, invariant under rescalings of the metric.
 
\ 
 
\noindent 
{\bf Definition.\ts}
{\it 
A noncompact complete Riemannian manifold $(M,g)$   
is called {\rm asymptotically flat} if $A(M)=0$. 
}
 
\  

Note that the 
mere condition of being asymptotically flat places in general
no restrictions whatsoever on the topology of a manifold.
For instance, by a result of Abresch (see [Ab])
{\sl any} noncompact surface carries
a complete and asymptotically flat Riemannian metric.

\
 
\noindent 
{\bf Definition.\ts} 
{\it 
A noncompact complete Riemannian manifold $(M,g)$ is said  
to have {\rm cone structure at infinity} if there is a  
metric cone $C$ with vertex $o$ such that the pointed Gromov--Hausdorff limit 
of $(M,\varepsilon{\cdot} g,p)$ exists
for any sequence of numbers $\varepsilon>0$ converging to zero
and such that this limit is isometric to $(C,o)$. 
}

\ 

As a Gromov--Hausdorff limit
of proper metric spaces,
that is, metric spaces such that any closed ball of finite radius is compact,
the cone $C$ which arises in the above definition
is in particular proper and locally compact.
 
\

Note that large classes of Riemannian manifolds have
cone structure at infinity.
In fact,
Kasue (see [K1]) showed that under certain lower curvature limitations 
(e.g., if for some $C<\infty$ and $\delta>0$ it holds that 
$K_x\ge -C/|px|^{2+\delta}$) 
a noncompact (complete) Riemannian manifold always has this property. 
Thus in particular
any noncompact Riemannian manifold
with faster-than-quadratic curvature decay
(that is, any noncompact Riemannian manifold for which 
there exists some $C$ and $\delta>0$ such that
$|K_x|\le C/|px|^{2+\delta}$),
and, especially, any noncompact Riemannian manifold
with nonnegative curvature
has cone structure at infinity.

On the other hand by Abresch's result one can easily construct
asymptotically flat surfaces $(S,g)$
such that the Gromov--Hausdorff limit 
of $(S,\varepsilon_n{\cdot} g,p)$ indeed depends on the choice of the 
sequence $\varepsilon_n\to 0$
and such that for some sequences this limit is not even a metric cone.
In particular, by considering products $\TT^{m-2}\times (S,g)$  
one thus obtains examples 
of asymptotically flat manifolds without cone structure at infinity 
in any dimension $m\ge 2$.
(Actually any such example we know of
looks on a big scale always two-dimensional;
for more on this see section $3$).

\

\noindent 
{\bf Theorem A.\ts}
{\sl 
Let $M$ be an asymptotically flat $m$-manifold.  Assume 
that $M$ has cone structure at infinity. 
Then  
 
\item{(i)}  
There exists an open ball $B_R(p)\i M$
such that $M\setminus B_R(p)$ is a disjoint 
union $\bigcup_i N_i$ of a finite number of ends;  
that is,  $N_i$ is a connected topological manifold with closed boundary 
$\partial N_i$
which is homeomorphic to $\partial N_i\times [0,\infty)$.
For each $N_i$ the limit $C_i=$
GH-$\lim_{\varepsilon\to 0} \varepsilon{\cdot} N_i$ exists.

\item{(ii)} If the end $N_i$ is simply connected, 
then $N_i$ is homeomorphic to $\SS^{m-1}\times [0,\infty)$.

In this case moreover the following holds:
 
\itemitem{(a)} if $m\not=4$, then $C_i$ is isometric to $\RR^m$;
 
\itemitem{(b)} if $m=4$, then 
$C_i$ is isometric to one of the following spaces:
$\RR^4$, $\RR^3$, or $\RR\times [0,\infty)$.
}

\

A finiteness of ends statement as in part (i) of Theorem $A$ 
was proved by Abresch (see [Ab]) in a related setting.

Part (ii)(a) of Theorem $A$
is new even in the special case of faster-than-quadratic curvature decay
(recall that, as noted above,
faster-than-quadratic curvature decay
implies cone structure at infinity).

Theorem $A$ generalizes here 
work of Greene, Petersen, and Zhu 
(see Theorem $1$ in [GPZ], 
where the same conclusion as in (ii)(a) was proved 
under the additional assumptions of 
faster-than-quadratic curvature decay
and nontriviality of the tangent bundle of $\partial N_i$.

\

When combined with the fact
that volume growth of exactly Euclidean order implies flatness,
Theorem $A$ also yields the following result,
which generalizes for manifolds of dimension $m\ne 4$ 
Theorem $2$ in [GPZ]
from faster-than-quadratic curvature decay to
asymptotical flatness with cone structure at infinity:

\

\noindent 
{\bf Corollary.\ts}
{\sl
Let $M$ be an asymptotically flat $m$-manifold of dimension $m\ne 4$
which has cone structure at infinity.
If $M$ has nonnegative Ricci curvature and one simply connected end,
then $M$ is isometric to $\RR^m$.
}

\

Part (ii)(b) of Theorem $A$ shows that 
four-dimensional manifolds play in Theorem $A$ a peculiar role
and indicates that in this dimension special phenomena can arise.
Indeed, Unnebrink (see [U]) showed that there are
examples of asymptotically flat $4$-manifolds
which have (cone structure at infinity and)  
a simply connected end $N_1$ 
such that, in the notation of Theorem $A$,  $C_1=\RR^3$.
It is not clear if there exists
an asymptotically flat $4$-manifold 
with cone structure at infinity
with a simply connected end $N_1$ so that $C_1=\RR\times [0,\infty)$.
(We actually conjecture that there is no such example;
see also section $3$.)

Note also that in dimension $2$ all ends are homeomorphic 
to $\SS^1\times [0,\infty)$, so
that the
ends of an asymptotically flat surface are never simply connected.

\

Combining Theorem $A$ and a result from [GP]
we obtain
proofs of statements of Gromov (see [BGS], p.59)
which till now have been treated in the literature 
(compare [D] and the references therein)
as conjectures,
and which completely classify 
the asymptotically flat manifolds  
with nonnegative sectional curvature:

\  
 
\noindent 
{\bf Theorem B.\ts}
{\sl  
Let $M$ be an asymptotically flat $m$-manifold with nonnegative sectional 
curvature.
Then the universal covering of $M$ is isometric to $\RR^{m-2}\times S$,  
where $S$ is an asymptotically flat surface. 
If, in particular, $M$ is simply connected and $m\ge 3$, 
then $M$ is isometric to $\RR^m$.
}

\

Assuming faster-than-quadratic curvature decay and 
assuming 
that the unit normal bundle of the soul of $M$ 
has nontrivial tangent bundle,
the second part of Theorem $B$ was proved by Drees ([D]).

As a direct consequence of Theorem $B$ one also 
obtains an affirmative answer
to a question of Hamilton
(see [H], \S 19; this paper also contains some nice relations
between asymptotical curvatures and the Ricci flow), 
which is equivalent to
the following one:

\

\noindent
{\sl 
Let $M$ be a complete noncompact Riemannian manifold
of dimension $m\ge 3$ with positive sectional curvature.
Is it true that $A(M)>0$ ?
}

\

That in odd dimensions the answer to this question is ``yes''
was already known and proved by Eschenburg, Schr\"oder, and Strake ([ESS]).

\ 
 

Our results obviously 
also have a relation to the positive mass conjecture;
in [GPZ] the reader will find explained the precise connections.

\ 

\

The main idea of the proof of Theorem $A$ can be described as follows:

\

Inside an end $N_i$ of an asymptotically flat manifold with cone
structure at infinity we construct a continuous family of ``spheres'',
which after rescaling have uniform curvature bounds and which
Gromov--Hausdorff converge
to the ``unit sphere'' in the cone $C_i$. 

To this continuous family of spheres we now apply two results from [PRT]:
The first says that any continuously collapsing family with bounded curvature 
contains an infinite ``stable''
subsequence. 
To this sequence then
a corollary of the Limit of Covering Geometry Theorem from [PRT] applies.

(This corollary actually also holds without using a stability assumption,
and the proof of Theorem $A$ can be given without
relying on [PRT], but instead on results from [PT],
see section $1$).

This in turn enables us 
to prove some inequalities for the ranks of certain homotopy
groups. These imply that in fact collapse is not possible except for 
the case where the dimension of the manifold is equal to $4$.
Therefore in all nice cases $C_i$ is nothing but $\RR^m$.

\

\

There is a vast amount of literature 
on noncompact complete Riemannian manifolds 
whose sectional curvature at infinity is zero
(and on many different specific 
ways in which the curvature is allowed to go to zero).
For a detailed account of what is known and wanted to be known 
about such spaces, the reader is recommended to look at 
the survey article [Gre] and the paper
[GPZ].
Here we just mention
(besides the references already given)
some papers in the field which are most closely related 
to the results of this note: [ES], [GW1], [GW2], [KS] and [LS].

\

The remaining parts of the paper are organized into a preliminaries,
a proof, and a problem section which contains further remarks
and several open questions.

\

We would like to thank Patrick Ghanaat 
for pointing out to us a simplified proof of
the sublemma in section $2$
as well as Luis Guijarro for useful comments.

\

\head 1. Preliminaries
\endhead 
 
\

In this section we review some results about manifolds which 
collapse with bounded curvature and diameter.
More on this can be found in the references given in [PRT] and [PT].

\

\noindent
{\bf Definition.\ts}
{\it
A sequence of metric spaces $M_i$ is called stable if 
there is a topological space $M$ and a sequence of metrics $d_i$ on $M$
such that
$(M,d_i)$ is isometric to $M_i$ and such that the metrics
$d_i$ converge as functions on $M\times M$ 
to a continuous pseudometric.
}

\

\noindent
{\bf Proposition (Continous Collapse implies Stability) ([PRT]).\ts} 
{\sl
Suppose that a simply connected manifold $M$  
admits a continuous  
one-parameter family of metrics $(g_t)_{0<t\le 1}$ with 
$\lambda\le
K_{g_t}\le 
\Lambda$ such that, as $t\to 0$, the family of metric spaces $M_t=(M,g_t)$ 
Hausdorff converges to a compact metric space $X$ of lower dimension.  
Then the  
family $M_t$ contains a stable subsequence. 
}

\

The version of the Limit of Covering Geometry Theorem from [PRT] 
we need in this paper
(it is straightforward to check 
that the proof given in [PRT] also proves the result below)
can be stated as follows:
 
\ 
 
\noindent
{\bf Theorem (Limit of Covering Geometry Theorem ([PRT]).\ts} 
{\sl 
Let $M_n$ be a stable
sequence of Riemannian $m$-manifolds with curvature bounds $|K(M_n)|\le1$ such 
that for $n\to\infty$ the sequence of metric spaces $M_n$ Hausdorff converges
to 
a compact metric space $X$ of lower dimension.   
Consider any sequence of points 
$p_n\in M_n$ and balls $B_n=B_{\pi/2}\i T_{p_n}$ which are equipped with the 
pull back metrics of the exponential maps 
$$\exp_{p_n}\colon T_{p_n}\to M_n.$$   
Assume that for any such converging subsequence $B_n\to B$,  
the limit $B$ has curvature $\ge0$ in 
the sense of Alexandrov. 

Then for any converging subsequence $B_n\to B$, the limit $B$  
has the same dimension as the manifolds $M_n$, 
and in a neighbourhood of its center,  
the metric on $B$ coincides with that of 
a metric product $\RR\times N$, where $N$ is a manifold with two-sided bounded 
curvature $0\le K(N)\le 1$ in the sense of Alexandrov. 
}

\ 

Our proof of Theorem $A$ will in fact only use the following corollary
of this theorem. At first sight this corollary looks almost obvious,
but it doesn't seem easy to adopt
any of the known proofs of injectivity radius estimates to this case.

\noindent
{\bf Corollary.\ts} 
{\sl Let $M_n$ be a (stable)
sequence of closed simply connected
Riemannian manifolds of dimension $m\ge2$ with curvature
$|K(M_n)|\le C$ and uniformly bounded diameters.   
Consider any sequence of points 
$p_n\in M_n$ and balls $B_n=B_{\pi/2\sqrt{C}}\i T_{p_n}$
which are equipped with the 
pull back metrics of the exponential maps 
$$\exp_{p_n}\colon T_{p_n}\to M_n.$$   
Assume that for any such converging subsequence $B_n\to B$,  
the limit $B$ has
at all interior points curvature $=1$.
Then the manifolds $M_n$ converge to a standard sphere.
}

\

The stability condition is actually not necessary for the above result
to hold.
This can be seen from the following
independent proof of the Corollary,
which does not use stability at all.
The proof itself is very short, but since it
uses the notion of
Grothendieck--Lipschitz convergence and Riemannian megafolds from [PT],
we decided to also incorporate the above [PRT] approach,
which might be easier to understand.

\

\noindent
{\bf Proof of the Corollary without stability assumption.}
 
The only nontrivial part is to establish a lower positive bound
for the injectivity radii of all manifolds $M_n$.

Since because of the Gauss--Bonnet theorem
the case $m=2$ is trivial, we may assume that $m\ge 3$.

If the manifolds $M_n$ would collapse, then, after passing to
a subsequence if necessary, one may assume
that
the manifolds $M_n$ Grothendieck--Lipschitz converge to a Riemannian megafold
${\frak M}$ which is not a manifold.

The assumptions of the Corollary imply that
the limit ${M}$ has constant curvature $=1$,
so that ${\frak M}=(\SS^m:G)$, where $G$ is a commutative group of isometries
of $\SS^m$.
However, by ([PT], Theorem A.7)
we have that $0\not=H^2_{dR}({\frak M})=H^2_{dR}(\SS^m)$, 
which is a contradiction.
$\qed$

\

\

\head 2. Proofs
\endhead

\ 
 
\noindent 
{\bf Proof of part (i) of Theorem A.}

The first statement of the theorem will follow from the fact  
that the distance function to $p$, dist$_p$,  
for sufficiently large values does not have any critical points.
 
By assumption, for any sequence of numbers $\varepsilon>0$ converging to zero, 
the pointed Gromov--Hausdorff limit of $(M,\varepsilon{\cdot} g,p)$ exists  
and is isometric to a locally compact metric cone $C$ with vertex $o$. 
The cone $C$ obviously has curvature $\ge 0$ (in the sense of Alexandrov)
everywhere except the origin $o\in C$.

Let us assume that there exists a sequence of points $x_n$  
such that $|px_n|\to\infty$ as $n\to\infty$, 
and such that each $x_n$ is a critical point for dist$_p$.  
 
Consider the sequence of rescaled manifolds $(M,g/|px_n|,p)$. 
By the assumption of the theorem, this sequence converges to $(C,o)$, 
and the points $x_n\in (M, g/|px_n|)$ (after passing to a subsequence)
converge to a point $x\in C$  
which has distance $1$ to the origin $o$. 
 
Since $C$ is a cone, we can consider $y:=2x\in C$. 
Choose a sequence of points $y_n\in (M,g/|px_n|)$ which converge to $y$,  
and consider minimal geodesics $x_ny_n$ from $x_n$ to $y_n$. 
Since $x_n$ is a critical point of dist$_p$, 
there is for each $n$ a minimal geodesic $px_n$ which makes an angle less than 
$\pi/2$ with the minimal geodesic $x_ny_n$. 
Therefore Toponogov's comparison theorem implies that 
$\lim_{n\to\infty}|py_n|/|px_n|\le \sqrt 2$. 
But obviously $\lim_{n\to\infty}|py_n|/|px_n|= |oy|/|ox|=2$, a contradiction. 

Thus for some $R>0$ the function
dist$_p$ does not have any critical points outside the open ball $B_R(p)$.
In particular, as follows from Morse theory for distance 
functions, see ([Grov], Cor.~1.9.), $M$ has finite topological type;
that is, $M$ is homeomorphic to the interior of a compact manifold with
boundary (which in our case is simply the closed ball $\bar B_R(p)$).

This also implies that the manifold $M$ has only finitely many ends.

Note that the cone $C_i$ is nothing but the closure of the connected component 
of
$C\backslash o$ that corresponds to $N_i$, 
in particular
for each $N_i$ the limit $C_i=$
GH-$\lim_{\varepsilon\to 0} \varepsilon{\cdot} N_i$ exists.

Thus part (i) of Theorem $A$ is proved. $\qed$

\ 
 
\noindent
{\bf Proof of part (ii) of Theorem A.}

The fact that $N$ is homeomorphic to $\SS^{m-1}\times [0,\infty)$ will follow
directly from the proofs of (ii)(a) and (ii)(b). 
Therefore we only need to prove these two statements.

Note that if dim\,$C_i=m=$ dim$M$, then parts (ii)(a) and (ii)(b) 
of the theorem are trivially true: 
 
Indeed, if so we have that the curvature of $C_i$ is zero everywhere except the 
origin.
We can assume that  $m\not=2$ (otherwise all ends would be homeomorphic to    
$\SS^1\times [0,\infty)$, and therefore in particular 
they would not be simply connected).
It 
follows that $C_i=\RR^m/F$,
where $F$ is a finite group of rotations which acts freely
on $\RR^m\backslash 0$.
Since $C_i \setminus B_1(o)$ is homeomorphic to $N_i$,
it follows that $F=\pi_1 (\partial N_i)$.
Since by assumption $\partial N_i$ is simply connected, $F$ must be trivial,
and thus for dim\,$C_i=m=$ dim$M$ our claims are proved,
since the above also implies that 
in this case $N_i$ is homeomorphic to 
$\SS^{m-1}\times [0,\infty)$.
 
\

From now on we will assume that dim $C_i<m$. 
 
\ 

We can view $C_i$ as a cone over its space of directions, 
$C_i=C(\Sigma_i)$, where $\Sigma_i$ is an Alexandrov space 
of curvature $\ge 1$ or dim$\Sigma_i=1$. $\Sigma_i$ can be viewed as a  
``unit sphere'' in $C_i$.  

\

We will first construct a continuous family of hypersurfaces
$S_{2-1,\varepsilon}$
in $(M,\varepsilon{\cdot} g)$ which collapse to 
$\Sigma_i$ such that the  
sectional curvatures of $S_{2-1,\varepsilon}$ stay uniformly bounded.  
The following construction is very close to one used by Kasue in [K2].
We will therefore only explain it here; all of its
details can be found in ([K2, $\S 2$]).
 
For each rescaling $(M, \varepsilon{\cdot} g,p)$, consider the sphere of radius $2$, 
$S_{2,\varepsilon}(p)\i (M, \varepsilon{\cdot}  g)$. 
Its principal curvatures for 
outcoming normal directions lie in the range $[-C(\varepsilon),\infty]$, 
where $C(\varepsilon)\to 1/4$ as $\varepsilon\to 0$. 
 
Next consider, for an inward direction (to $p$)  
the equidistant hypersurface $S_{2-1,\varepsilon}$
at distance $1$ to $S_{2,\varepsilon}(p)$. 
Then $S_{2-1,\varepsilon}$ has 
uniformly bounded principal curvatures which in fact lie in the range 
$[-C'(\varepsilon),C'(\varepsilon)]$, where $C'(\varepsilon)\to 1$ as $\varepsilon\to 0$.  

Therefore, since $M$ is asymptotically flat,
 $S_{2-1,\varepsilon}$ has uniformly bounded sectional curvature as
$\varepsilon\to0$. 
For sufficiently small $\varepsilon$ it follows that $S_{2-1,\varepsilon}$
(equipped with the induced intrinsic metric) 
is a continuous family which, as $\varepsilon\to 0$, 
collapses to~$\Sigma_i$. 

\ 

\noindent
{\bf Key Lemma.\ts} 
{\sl
Take any sequence of points
$p_{\eps_n}\in S_{2-1,\eps_n},\  \eps_n\to 0$. 
Consider the balls $B_n=B_{1}(0)\i T_{p_{\eps_n}}(S_{2-1,\eps_n})$,
equipped with
the pull back metrics.

Then as $n\to \infty$, the $B_n$ Lipschitz converge to the 
ball of radius $1$ in 
$\SS^{l-1}\times\RR^{m-l}$, for some fixed $l$ depending on $M$.

Moreover $\Sigma_i=\widetilde \SS^{l-1}/A$, where  $A$ is an Abelian group 
of isometries of $\widetilde \SS^{l-1}$ 
(here by $\widetilde \SS^{l-1}$ we understand
the standard $l-1$-sphere if $l\ge 3$, $\RR$ if $l=2$, 
and a point if $l=1$).
}

\ 

The proof of the Key Lemma will be given below. 
Let us now continue with the proof of part (ii) of Theorem $A$:

\ 
 
Obviously all $S_{2-1,\varepsilon}$ are homeomorphic to $\partial N_i$ and 
therefore simply connected.
Now applying the Corollary in section $1$
we see that $l<m$. 

Using that $\partial N_i$ is simply connected we can moreover show that 
the group $A$ in the Key Lemma is connected:
Let $A^o$ be the identity component of $A$.
Then
$\widetilde \Sigma_i=\SS^{l-1}/A^o$ is a branched covering of
$\Sigma_i$, 
and it is easy to see that one can find
a covering $\widetilde {\partial N_i}\to \partial N_i$ which is a lifting of
$\widetilde \Sigma_i\to\Sigma_i$. But since $\partial N_i$ is simply connected
we have that $\widetilde \Sigma_i=\Sigma_i$. 
Therefore $A^o=A$; that is,
$A$ is connected.

\

Now if $l=1$, then $\Sigma_i$ is a point, so that $\partial N_i$ must be an 
infranil manifold. But any infranil manifold  
has infinite fundamental group, which contradicts the fact that 
$\partial N_i$ is simply connected.

If $l=2$ it follows that $\Sigma_i$ is 
homeomorphic to a point or $\RR$, since $A$ is connected.
The first case cannot occur by the above reasoning, and the second 
contradicts that $C$ is locally compact.

\ 

Therefore the only serious case to deal with is the case $l\ge 3$.

From the above we have that in this case
$\Sigma_i$ is isometric to $\SS^{l-1}/\TT^{k'}$.

Since for $\eps\to 0$ the hypersurfaces $S_{2-1,\eps}$ collapse
to $\Sigma_i$ and since $S_{2-1,\eps}$ is homeomorphic to $\partial N_i$,
we know that  $\Sigma_i$ is homeomorphic to $\partial N_i/\TT^k$ and 
that this homeomorphism can be chosen to preserve the natural stratifications
of these spaces.

\

Let us now do some topological calculations:

Let $O_{\TT^k}$ be a regular orbit of the $\TT^k$ action on $\partial N_i$.
Consider the relative homotopy sequence of pairs
$$ \pi_2(\partial N_i,O_{\TT^k})\to \pi_1(\TT^k)=\ZZ^k\to \pi_1(\partial N_i)=0.$$ 
Therefore rk$_\QQ\,\pi_2(N_i,O_{\TT^k})\ge k$. 

\
 
Next consider the corresponding homotopy sequence for $\SS^{l-1}$:
$$0=\pi_2(\TT^{k'})\to \pi_2(S^{l-1})\to \pi_2(\SS^{l-1},O_{\TT^{k'}})\to 
\pi_1(\TT^{k'})=\ZZ^{k'}\to 
\pi_1(\SS^{l-1})$$ 
Therefore rk$_\QQ\,\pi_2(\SS^{l-1},O_{\TT^{k'}})=k'+$ 
rk$_\QQ\,\pi_2(\SS^{l-1})$. 
 
On the other hand one has that $k=$ dim$\partial N_i-$ dim$\Sigma_i$ and 
$k'=l-1-$ dim$\Sigma_i$. 

\

\

Let $\Sigma_i^\#$ denote $\Sigma_i$ with the singular sets removed.
Consider now the following three cases:

\

1. $\Sigma_i$ has no boundary. 
Then obviously
$$\hbox{rk}_\QQ\,\pi_2(\SS^{l-1},O_{\TT^{k'}})
=\hbox{rk}_\QQ\,\pi_2( \Sigma_i^\#)
=\hbox{rk}_\QQ\,\pi_2(\partial N_i,O_{\TT^k}).$$

\

2. $\partial \Sigma_i$ has one component. 
Then 
$$\hbox{rk}_\QQ\pi_2(\SS^{l-1},O_{\TT^{k'}})=\hbox{rk}_\QQ\pi_2( \Sigma_i^\#)+1=\hbox{rk}_\QQ\pi_2(\partial N_i,O_{\TT^k}).$$

\

3. $\partial \Sigma_i$ has more than one component.

\ 

In both case $1$ and case $2$ it follows that 
$k'\ge k-$ rk$_\QQ\pi_2(\SS^{l-1})$, and hence  
$$m\le l+\text{rk}_\QQ\pi_2(\SS^{l-1}).$$

However, this contradicts the fact that $l<m$, 
except if $m=4, l=3$.
In this particular case it follows that $\Sigma_i=\SS^2/A^{rot}$.
Therefore,
since $\Sigma_i$ has not more than one boundary component,
we have that
$A^{rot}$ is trivial and $\Sigma_i=\SS^2$.
Thus $C_i=C(\Sigma_i)=\RR^3$
(and that this indeed can happen was shown in [U]).

\

Case $3$ can only occur
if $\Sigma_i$ is homeomorphic to $[0,1]$.
Then, since $N_i$ is simply connected, it must hold that $k$, $k'\le 2$.
Since  the $\TT^k$ action on $\partial N_i$ has empty fixed point set,
we have that $k=2$,
and since $l<m$, we have that $k'=1$. 
Therefore $m=4$, $l=3$ and $\Sigma_i$ is isometric
to $\SS^2/\SS^1=[0,\pi]$ , so that $C_i=C(\Sigma_i)=\RR\times [0,\infty)$.  

The proof of Theorem $A$ is complete. \qed

\ 

\

\noindent
{\bf Proof of Theorem B.}

Let $M$ be an asymptotically flat $m$-manifold with nonnegative sectional
curvature. Then $M$ has cone structure at infinity,
and by [GP] the soul $S$ of $M$ is flat.
This forces the universal cover $\tilde M$ of $M$
to split isometrically as a Euclidean part, coming from the soul $S$,
and a nonnegatively curved complete manifold $F$ homeomorphic to $\RR^k$.

Now $F$ is also asymptotically flat
and has one end $\SS^{k-1}\times [0,\infty)$.
Therefore by Theorem $A$, if $k\ne 2,4$, then the cone at infinity of $F$
is isometric to $\RR^k$. 
Since by Toponogov's Comparison Theorem
any line in the cone corresponds to a line
in $F$, it follows from the Toponogov Splitting Theorem
that $F$ itself is isometric to $\RR^k$.

Thus to finish the proof we must only exclude the case $k=4$.
By Theorem $A$, if $k=4$ we have that $C=$ 
GH-$\lim_{\varepsilon\to 0} \varepsilon{\cdot} F$ is isometric to one of the following:
$\RR^4$, $\RR^3$ or $\RR\times [0,\infty)$.
In all of these cases we have
that $C$ contains a line, and therefore $F$ splits isometrically as 
$\RR\times F'$.
But since $F$ is asymptotically flat it follows that $F'$ is flat,
and therefore $F$ is isometric to $\RR^4$. $\qed$

\

\noindent
{\bf Proof of the Key Lemma.} 

Consider a $\nu$-neighbourhood $U\supset \Sigma_i\i C_i$.  
From the results of [CFG] (see section $1$ of [PRT], where also further
references can be found)
we have an $N$-structure 
$\pi\colon  E_\varepsilon\to U$, where  
$E_\varepsilon$ is a subset of $(M, \varepsilon{\cdot} g)$
containing the hypersurface $S_{2-1,\varepsilon}$.  
Since $E_\varepsilon$ is homotopically equivalent to $\partial N_i$,  
it follows that $E_\varepsilon$ is simply connected.
Therefore the $N$-structure is given by
an almost isometric smooth $\TT^k$-action without fixed points 
(see again section $1$ in [PRT]).

Now take a 
point $x\in \Sigma_i\i C_i$ (so $|ox|=1$) and consider a spherical  
neighbourhood of $U_x\ni x$.  
Consider the preimage $V_\varepsilon=\pi^{-1}(U_x)\i E_\varepsilon$ and let   
$\widetilde V_\varepsilon$ be its universal Riemannian covering. 
Then the $\TT^k$-action induces an almost isometric $\RR^k\times F$ action on  
$\widetilde V_\varepsilon$, where $F$ is a finite 
Abelian group.  

From [CFG] one has a uniform bound for  
the injectivity radius of $\widetilde V_\varepsilon$, so that,  
as $\varepsilon\to 0$, 
$\widetilde V_\varepsilon$ converges to a flat manifold $\widetilde V_0$ 
with boundary and isometric $\RR^k\times F$ action 
(for the convergence claim see the first part of Lemma 2.1.4 in [PRT]).

Since the interior of $\widetilde V_0$ is flat,  
there exists a map $\widetilde V_0\to \RR^m$ 
which is for all interior points a local isometry.  
Therefore the $\RR^k\times F$ action on 
$\widetilde V_0$ can be extended to an action of whole $\RR^m$, 
and the local factors $U\i \RR^m/\RR^k$ are isometric to local 
branched coverings of subsets of $C_i$.  
(Here by local  
factors we understand factorizing $U$ by the connected components  
of the $\RR^k$-orbits  
in $U$, as is illustrated in the following picture). 
 
\  

\
 
\input epsf 
\vskip0.0in 
\epsfxsize2.0in 
\hfil\epsfbox{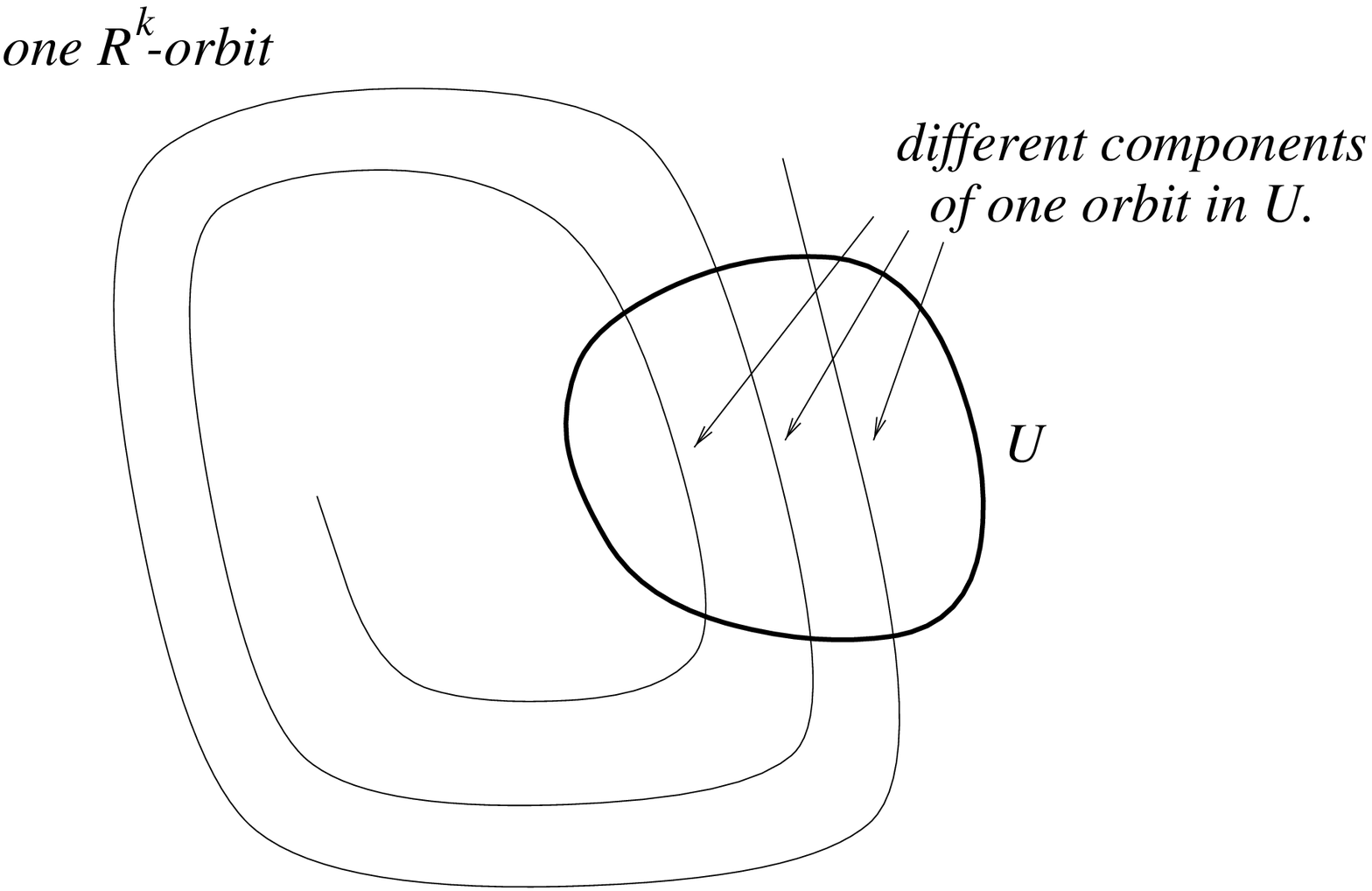}\hfil 
 
\ 

\

Now the above group $\RR^k$ can be regarded as an Abelian group of isometries  
of Euclidean space $\RR^m$.
We will show that in our case $\RR^k$ actually
splits into a direct sum of translations and rotations.

\

To this means first note the following:

\
 
\noindent
{\bf Sublemma.\ts} 
{\sl Let a connected Abelian group $H$ act on Euclidean space
$\RR^m$ by isometries.  
Then one can represent $\RR^m$ as an orthonormal sum $V\oplus W$ such that  
$H$ is contained in a direct sum of translations and rotations,
$$H<A^{tr}\oplus A^{rot},$$
so that the following holds:
The group $A^{tr}=V$ is the group 
consisting of all parallel translations along $V$, 
and $A^{rot}\i O(W)$ is an Abelian subgroup of rotations of $W$. 
}

\

\noindent
{\bf Proof of the Sublemma.}

By [Al] one orbit of $H$ is an affine subspace $V$
(in fact, such an orbit corresponds to the origin $o$ of $C_i$).
Choose the origin of 
affine space $\RR^m$ so that it is contained in this subspace.
Each element $\alpha\in H$ can be viewed as 
$(r_\alpha,\phi_\alpha)\in V\times O(m)$, such that 
$\alpha(x)=r_\alpha+\phi_\alpha(x)$ for any
$x\in \RR^m$.

Then $V$ can be viewed as the set of all pure translations of $H$,
$$A^{tr}=V=\{r:(r,\phi)\in H\ \hbox{for some}\  \phi\in O(m)\}.$$ 
Let $A^{rot}:=\{\phi:(r,\phi)\in H$ for some $r\}$ be the group of pure 
rotations of $H$.
If each $\phi\in A^{rot}$ acts trivially on $V$,
then obviously $H<A^{tr}\oplus A^{rot}$, which is exactly what we want.

Therefore we only have to prove that for any $\phi\in A^{rot}$ and any $v\in V$
we have that $\phi(v)=v$.

Take any $(r,\phi)\in H$ and $v\in V$. 
For all $n\in\NN$ there exists $\phi_n\in O(m)$ such that 
$(nv,\phi_n)\in H$.
Since $H$ is Abelian, it follows that 
$(r,\phi)(nv,\phi_n)= (nv,\phi_n)(r,\phi)$ 
and therefore $nv+\phi_n r = r+\phi nv$. 
Dividing by $n$ and letting $n\to\infty$ thus implies
$\phi v = v$. $\qed$

\

Thus our group $\RR^k$ is contained in a direct sum 
$A^{tr}\oplus \tilde A^{rot}$, where $\tilde A^{rot}$ is universal covering 
Lie group of 
$A^{rot}$.
Now note that since the local factors by $\RR^k$ have a cone structure, 
$\RR^k$ moreover itself splits as $\RR^k=A^{tr}\oplus \tilde A^{rot}$: 

Indeed,
since the local factors $U/\RR^k$
admit a cone structure, in radial directions 
their sectional curvatures must be zero. 
But this is impossible unless $\RR^k$ is itself a direct product 
$\RR^k= \tilde H= A^{tr}\oplus \tilde A^{rot}$:

To prove this, 
we only have to show that 
(in the notation of the Sublemma) it holds that $A^{rot}\i H$. 
Assume that this is wrong.
Then we can find a ray $c\colon [0,\infty)\to \RR^m$ 
which is orthogonal to $V$, and there will be an element 
$\alpha\in a^{rot}$ in the Lie algebra of $A^{rot}$
which is not contained in the Lie algebra of $H$,
so that $\alpha$ defines a linear Jacobi field on the ray $c$
which can assumed to be non-zero.

Consider now the projection $\bar c$ of $c$ along some local factor. 
Then $\bar c$ is a piece of a ray in the cone $C$ and the
projection $\bar J$ of 
the field $J$ is also a Jacobi field. 
But since $C$ is a cone, any 
Jacoby field along $\bar c$ must be linear.
On the other hand it is straightforward to show that $|\bar J(t)|$ is a  
strictly concave function, and this is a contradiction.

\

Therefore the   
local factors $W/A^{rot}$ are isometric to local branched coverings
of $C_i$ (everywhere except the origin).
Thus
$C_i\setminus o$ is isometric to a factor of its universal covering,
$\widetilde{(W\setminus 0)}/A$, by an Abelian Lie group $A$.
Restricting this last isometry to the unit spheres of both cones 
it follows that $\Sigma_i=\tilde \SS^{l-1}/A$, and 
the second part of the Key Lemma is proved.

\

Let $\rho\colon \widetilde V_\varepsilon\to V_\varepsilon$ be the covering map and
$\widetilde S_{2-1,\varepsilon}=\rho^{-1}(S_{2-1,\varepsilon})$.
It converges to the preimage of $\Sigma_i$  
under the map $V_0\to V_0/A=U_x\i C_i$, so that it locally coincides with
the cylinder  
$V\times \SS^{l-1}$, where $\SS^{l-1}$ is the unit sphere in $W$.   
Therefore, since $x\in\Sigma_i$ can be chosen to be arbitrary, the covering 
geometry of $S_{2-1,\varepsilon}$ converges to the one of $V\times \SS^{l-1}$, 
and this finishes the proof of the first part of the Key Lemma.  $\qed$

\

\

\head 3. Remarks and Open Questions
\endhead 

\ 

\noindent
{\bf Question 1.\ts} Let $M$ be an asymptotically flat manifold, 
and let the sequence $(M,\varepsilon_n{\cdot}g,p)$ converge to $(G,o)$ 
as $\varepsilon_n\to 0$.
Assume that dim$G\ge 3$ and that $G\backslash o$ 
has only one connected component.

Is it true that $G$ is a metric cone with origin $o$?

\ 

A positive answer to this question could possibly
lead to a general classification 
of asymptotically flat manifold of higher dimension.
To obtain such a classification
is particularly interesting because of the fact
that
Gromov (see [Grom], p.96 and also [LS]) 
showed that any (smooth paracompact) noncompact manifold
$M$ admits a complete Riemannian metric whose asymptotic
curvature satisfies $A(M)\le C$,
where $C$ depends only on the dimension of $M$.

\

\noindent
{\bf Question 2.\ts} 
Does there exist in each dimension $m$
a positive constant $C(m)$ such that any noncompact
complete Riemannian $m$-manifold $M$ with $A(M)\le C(m)$
admits an asymptotically flat metric?

\ 

Note that the answer (positive or negative) 
to the following question would give
a complete classification
of the cone structures at infinity of simply connected ends
of asymptotically flat manifolds:

\ 

\noindent
{\bf Question 3.\ts} Can the cone $\RR\times [0,\infty)$ 
be a cone at infinity of a simply connected end of an
asymptotically flat $4$-manifold?

\ 

It seems not possible to obtain such an 
example by a direct generalization of Unnebrink's example. 
Namely, one can exchange 
the Berger spheres 
$\SS^3_{f(t),h(t)}$ (in the notation of [U]) 
in Unnebrink's example by $\SS^3_{a(t),f(t),h(t)}$, 
where the number $a(t)$ describes along which one-dimensional subgroup 
of the $\TT^2$-action 
on the standard $\SS^3$ we shrink the distance 
(so $\SS^3_{a,f,h}$ is a Berger sphere if $a=\pm1$).
But direct calculation then shows that there is no triple of functions
$a,f,h$ which
would give an asymptotically flat 4-manifold with $\RR\times [0,\infty)$
as a cone at infinity.

\

However, on the other hand, if one would take as $a$ a constant 
which is close to $1$, then
as a result one obtains an end $N$ whose asymptotic
curvature $A(N)$ is arbitrarily small and 
which has $\RR\times [0,\infty)$
as cone at infinity.

\ 

\

\noindent
{\bf Remark.\ts} 
The same arguments as the one which we used in the proof of Theorem $A$
actually also make it possible to characterize the 
cones at infinity of
complete noncompact manifolds whose asymptotic curvature is {\sl small}:

Namely, if for some given sequence of simply connected 
$m$-dimensional ends $N_n$ with $A(N_n)\to 0$ as $n\to\infty$
their cones at infinity $C_n$ Gromov--Hausdorff converge to some metric space
$C$ (which then must be a cone), 
then for 
sufficently large $n$ it holds that 
$N_n$ is homeomorphic to $\SS^{m-1}\times [0,\infty)$ and moreover the following
is true:
 
\itemitem{(a)} if $m\not=4$, then $C$ is isometric to $\RR^m$;
 
\itemitem{(b)} if $m=4$, then 
$C$ is isometric to one of the following spaces:
$\RR^4$, $\RR^3$, or $\RR\times [0,\infty)$.

\ 

The above modification of the Unnebrink construction for constant $a$
shows that for manifolds with small asymptotic curvature
all cones which are mentioned in part (b) 
actually do arise. 

\

As a last point we would like to mention 
that the methods we used in this paper do not distinguish
between spaces which are asymtotically flat and sequences of spaces
whose asymptotic curvature goes to zero. 

Therefore, no matter 
how special our question 
whether 
$\RR\times [0,\infty)$ 
can be a cone at infinity of a simply connected end of an
asymptotically flat $4$-manifold
might at first sight look like,
any negative answer to it 
will require more sensitive collapsing techniques.

\

\

\

\

\

\enddocument